\documentclass[reqno,11pt]{amsart}
\usepackage{amssymb, amsthm, amsfonts, amsgen, stmaryrd}
\usepackage{slashed}
\usepackage[english]{babel}
\usepackage{fullpage}
\usepackage[centertags]{amsmath}
\usepackage{indentfirst}
\usepackage{fancyhdr}
\usepackage[dvips]{graphicx}
\usepackage{psfrag}
\usepackage{enumerate}
\numberwithin{equation}{section}
\usepackage[latin1]{inputenc}

\usepackage{color}

\title[Hodge - De Rham   operator on $\SU$] {A Hodge - de Rham  Dirac operator on the quantum   ${\rm SU}(2)$}
\date{18 september 2016}
\author{Fabio Di Cosmo}
\author{Giuseppe Marmo}
\author{Juan Manuel P\'erez - Pardo}
\author{Alessandro Zampini}
\address{Dipartimento di Scienze Fisiche, Universit\`a di Napoli Federico II and INFN - Sezione di Napoli, Via Cintia - 80126 Napoli, Italy} 
\address{Departamento de Matem\'aticas, Universidad Carlo III de Madrid, Avda. de la Universidad 30, 28911 Legan\'es, Madrid}
\email{dicosmo@fisica.unina.it}
\email{marmo@na.infn.it}
\email{jmanuelperpar@gmail.com} 
\address{ Mathematics Research Unit, University of Luxembourg, 6 Rue Richard Coudenhove - Kalergi - 1359  Cit\`e de Luxembourg - Luxembourg}
\email{azampini@gmail.com}

\newcommand{\nn}{\nonumber}

\newcommand{\dd}{{\rm d}}
\newcommand{\ca}{\mathcal{A}}

\newcommand{\ch}{\mathcal{H}}

\newcommand{\cq}{\mathcal{Q}}

\newcommand{\cu}{\mathcal{U}}        
\newcommand{\SU}{\mathrm{SU}_q(2)}  
\newcommand{\ASU}{\ca(\mathrm{SU}_q(2))}  
\newcommand{\su}{\cu_q(\mathfrak{su}(2))}  

\newcommand{\eps}{\varepsilon}      
\newcommand{\hs}[2]{\left\langle #1,#2\right\rangle}  

\newcommand{\lt}{{\triangleright}}    
\newcommand{\rt}{{\triangleleft}}
\newcommand{\IC}{{\mathbb C}} 
\newcommand{\IR}{{\mathbb R}} 

\newcommand{\figureheight}{8cm}
\newcommand{\putfig}[2]{\begin{figure}[htp]
        \special{isoscale c:/itex/texfig/#1.wmf, \the\hsize \figureheight}
        \vspace{\figureheight}
        \caption{#2}\label{fig:#1}
        \end{figure}}
\newcommand{\pictureheight}{4cm}
\newcommand{\putpicture}[2]{\begin{figure}[htp]
        \special{isoscale c:/itex/texfig/#1.wmf, \the\hsize \pictureheight}
        \vspace{\pictureheight}
        \caption{#2}\label{fig:#1}
        \end{figure}}

\newcommand{\beqa}{\begin{eqnarray}}
\newcommand{\eeqa}{\end{eqnarray}}
\newcommand{\beq}{\begin{equation}}
\newcommand{\eeq}{\end{equation}}

\newcommand{\del}{\partial}

\newcommand{\oz}{\omega_{z}}
\newcommand{\op}{\omega_{+}}
\newcommand{\om}{\omega_{-}}
\newcommand{\ot}{\otimes}

\newcommand{\gm}{\mathrm{g}}
\newcommand{\cx}{\mathcal{X}}

\newcommand{\R}{{\mathbb{R}}}
\newcommand{\Z}{{\mathbb{Z}}}
\newcommand{\N}{{\mathbb{N}}}
\newcommand{\C}{\mathbb{C}}

\newcommand{\Cl}{\mathrm{Cl}}

\newcommand{\D}{\mathcal{D}}
\newcommand{\bu}{\bar{u}}
\newcommand{\bv}{\bar{v}}

\newcommand{\A}{\mathcal{A}}

\begin{document}

\thispagestyle{empty}

\begin{abstract}
We describe how it is possible to describe irreducible actions of the Hodge - de Rham Dirac operator upon the exterior algebra over the quantum spheres $\SU$ equipped with a three dimensional left covariant calculus.

\end{abstract}


\maketitle
\tableofcontents


\section{Introduction: the Hodge - de Rham operator on a manifold}
\label{sec:hd}

Let $M$ be a finite $N$-dimensional smooth manifold.  We denote by  $\Lambda(M)\,=\,\oplus_{k=0}^N\Lambda^k(M)$ the exterior algebra over it,  equipped with the Cartan exterior calculus  $(\Lambda(M), \wedge, \dd, i_X, L_X)$ with respect to the wedge product, where $\dd:\Lambda^k(M)\to\Lambda^{k+1}(M)$ is the exterior differential, $i_X:\Lambda^{k}(M)\to\Lambda^{k-1}(M)$ the contraction along the vector field $X$  defined on  $M$ and $L_{X}:\Lambda^k(M)\to\Lambda^k(M)$ the Lie derivative along $X$.  
On a local chart domain with a coordinate system $\{x^i\}_{i=1,\ldots,N}$  an exterior $k$-form is written as $\Lambda^k(M)\,\ni\,\phi\,=\,(1/k!)\phi_{a_1\cdots a_k}\dd x^{a_1}\wedge\ldots\wedge \dd x^{a_k}$ with $\phi_{a_1\cdots a_k}\,\in\,\Lambda^0(M)\sim\mathcal{F}(M)$. We also denote by $\mathcal{X}(M)$ the set of vector fields over $M$, with $\mathcal{M}\ni X=X^a\del_a$ along $\{\del_a\}_{1,\ldots,N}$, the local vector fields basis dual to $\{\dd x^a\}_{a=1,\ldots,N}$, with $X^a\,\in\,\mathcal{F}(M)\,=\,{\mathcal A}$.    

Let $M$ be equipped with a metric tensor $g$, whose local expression is  $g=g_{ab}\dd x^a\otimes\dd x^b$, or equivalently $g\,=\,g^{ab}\del_a\otimes\del_b$ with $g^{ab}g_{bc}=\delta^a_b$: setting\footnote{The reader may refer to  the well known monographs \cite{fof, fried, lami, nicolae}    for a more general description on spin structures and Dirac operators on a manifold. The book \cite{gbvf} and the lecture notes \cite{landsman,varillydosg}  describe an algebraic setting for the study of Dirac operators in both  classical and  non commutative geometry.}
\beq
\label{clm}
 \phi\vee\phi^{\prime}\,=\,\sum_s\,\frac{(-1)^{\tiny{\left(\begin{array}{c}s \\ 2 \end{array}\right)}}}{s!}g^{{a_1}{b_1}}\,\cdots\,g^{{a_s}{b_s}}(\check{\gamma}^s\{i_{a_1}\,\cdots \,i_{a_s}\,\phi\})\wedge\{i_{b_1}\,\cdots\,i_{b_s}\,\phi^{\prime}\}, 
\eeq
where $\phi,\,\phi^{\prime}$ are elements in $\Lambda(M)$ and $\check{\gamma}(\phi)\,=\,(-1)^k\phi$ for $\phi\,\in\,\Lambda^k(M)$, amounts to realise on $\Lambda(M)$ the well known Clifford product corresponding to the given metric $g$. It is immediate to see that one has
\begin{align}
\label{defcli}
\dd x^a\vee\dd x^b+\dd x^b\vee \dd x^a&=2 g^{ab}, \\ 
\dd x^a\vee\dd x^b&=\dd x^a\wedge\dd x^b+g^{ab}. \label{defcl}
\end{align}
The set $(\Lambda(M), \wedge, \vee, i_X, \dd)$ is usually called the K\"ahler - Atiyah algebra over $(M,g)$. The unital algebraic structure $(\Lambda(M), \vee)$ 
is the Clifford algebra $\Cl(M, g)$ on $M$ corresponding to the metric tensor $g$: it gives, after K\"ahler, an \emph{inner} product over an exterior algebra.   Notice that the Clifford algebra $\Cl(M, g)$ acts upon $\Lambda(M)$, with such an action $\Phi\,:\,\Lambda(M)\,\to\,{\rm End}(\Lambda(M))$   being generated by
\beq
\Phi(\dd x^a)\,:\,\phi\,\mapsto\,\dd x^a\wedge\phi\,+\,g^{ab}i_{\del_b}\phi.
\label{claac}
\eeq
The Levi Civita connection corresponding to 
$(M,g)$ is defined on $\mathfrak{X}(M)$  and dually on  $\Lambda^1(M)$ via 
\begin{align}
&\nabla\,:\,\Lambda^1(M)\,\to\,\Lambda^1(M)\,\otimes_{\A}\,\Lambda^1(M), \nn \\
&\nabla\,:\,\mathcal{X}(M)\,\to\,\Lambda^1(M)\,\otimes_{\A}\,\mathcal{X}(M),
\end{align} 
with
\begin{align}
\nabla(\dd x^a)&=\,-\dd x^s\,\otimes\,\Gamma^a_{sb}\dd x^b, \nn \\
\nabla(\del_a)&=\,\dd x^s\,\otimes\,\Gamma^b_{sa}\del_b
\label{nabd}
\end{align}
where $\Gamma^a_{bc}$ are the Christoffel symbols of the 
connection,   
\beq
\label{Cs}
\Gamma_{ji}^m\,=\,\Gamma_{ij}^m\,=\,\frac{1}{2}g^{mk}\left(\del_jg_{ki}\,+\,\del_ig_{kj}\,-\,\del_kg_{ij}\right).
\eeq
The action of $\nabla$ is extended to $\Lambda(M)$  requiring it to satisfy the Leibniz rule with respect to the wedge product. For any  $\phi\,\in\,\Lambda^k(M)$, its exterior covariant derivative is defined \cite{ka} via
\beq
\label{eD}
\mathfrak{D}\phi\,=\,\dd x^a\,\wedge\,\nabla_{a}\phi,
\eeq
where $\nabla_a\,=\,\nabla_{\del_a}$. The action of the K\"ahler  (also called Hodge - de Rham) operator is defined as:
\beq
\label{did}
\mathcal{D}\phi\,=\,\dd x^a\,\vee\,\nabla_a\phi.
\eeq 
From \eqref{clm} one immediately sees that
\beq
\label{dve}
\D\phi\,=\,\mathfrak{D}\phi\,+\,g^{ab}i_b\nabla_a\phi.
\eeq
A direct computation shows that $\mathfrak{D}\phi\,=\,\dd\phi$, so one has
\beq
\label{die}
\D\phi\,=\,\dd\phi\,+\,g^{ab}i_b\nabla_a\phi:
\eeq
if $\phi\,\in\,\Lambda^k(M)$, the image $\D\phi$ has a component of degree $k+1$ and a component of degree $k-1$. Notice that $\D\phi$ has a well defined parity with respect to the $\Z_2$ grading of $\Lambda(M)$. One can prove that, for $\phi\,\in\,\Lambda^k(M)$, 
\beq
\label{difi}
\D\phi=\dd\phi\,+\,(-1)^{N(k-1)}\star\dd\star\phi
\eeq
where one has introduced the Hodge duality operator corresponding to the metric $g$, i.e. an $\A=\mathcal{F}(M)$-bimodule map $\star\,:\,\Lambda^k(M)\,\to\,\Lambda^{N-k}(M)$ defined on a basis by
\beq
\label{hod}
\star(\dd x^{a_1}\wedge\cdots\wedge\dd x^{a_k})\,=\,g^{a_1b_1}\,\cdots\,g^{a_kb_k}i_{b_k}\,\cdots\,i_{b_{1}}\,\mu
\eeq
where $\mu\,=\,|g|^{1/2}\dd x^1\wedge\,\cdots\,\wedge\,\dd x^N$ is the invariant volume form 
corresponding to the metric $g$ (here $|g|$ denotes the determinant of the matrix $g_{ab}$).  The  relation 
\beq
\label{hos}
\star^2(\phi)\,=\,(-1)^{k(N-k)}({\rm sgn}\,g)\phi
\eeq
holds, for $\phi\in\Lambda^k(M)$, with $({\rm sgn}\,g)$ the signature of the metric tensor $g$. For riemannian (i.e. positive definite metrics $g$), the previous relation allows to write the action \eqref{difi} as
\beq
\label{kmtco}
\D\phi=\dd\phi-(-1)^k\star^{-1}\dd\star\phi.
\eeq
Via a metric $g$ on a manifold, one can introduce a scalar product among exterior forms. For homogeneous $\phi$ and $\phi^{\prime}$ with the same degree we set 
\beq
\label{scp}
\langle\phi\mid\phi^{\prime}\rangle\,=\,\int_{M}\phi\,\wedge\,\star\phi^{\prime},
\eeq
where the integration measure on $M$ comes from the volume form $\mu$ previously defined. 
Homogeneous  forms $\phi, \phi^{\prime}$ are defined orthogonal if their degrees are different (this amounts to say that $\Lambda^k(M)$ is orthogonal to $\Lambda^{k^{\prime}}(M)$ for $k\,\neq\,k^{\prime}$). 
Given any $\alpha\,\in\,\Lambda^k(M)$ and $\beta\,\in\,\Lambda^{k+1}(M)$, upon defining 
$\dd^*\,=\,(-1)^{N(k-1)+1}\star\dd\star\,:\,\Lambda^{k}(M)\,\to\,\Lambda^{k-1}(M)$ one proves that
\beq
\label{agd}
\langle\dd\alpha\mid\beta\rangle\,=\,\langle\alpha\mid\dd^*\beta\rangle\,+\,\int_M\dd(\alpha\wedge\star\beta).
\eeq
This relation shows that the operator $\dd^*$ can then be considered -- under suitable conditions on its domain --  the adjoint of $\dd$ in $\Lambda(M)$ with respect to the scalar product \eqref{scp}. The action \eqref{difi} of the Dirac operator
can be cast in the form
\beq
\label{diag}
\D\phi\,=\,\dd\phi\,-\,\dd^*\phi:
\eeq
from \eqref{hos} we have that $\dd^{*}\dd^*\,=\,0$, making it   immediate to prove the following relation
\beq
\label{laD}
\D^2\phi\,=\,-(\dd\dd^*+\dd^*\dd)\phi\,=\,(-1)^{kN}\{\star\dd\star\dd\phi\,+\,(-1)^N\dd\star\dd\star\phi\}
\eeq
for $\phi\,\in\,\Lambda^k(M)$. The square of the Dirac operator gives a version of the Hodge - de Rham Laplacian on $\Lambda(M)$.

It is easy to see that the action \eqref{claac} of the Clifford algebra $\Cl(M, g)$ upon $\Lambda(M)$   is highly reducible. 
Since the Clifford algebra $(\Lambda(M), \vee)$ is semi simple (simple for even $N$), its finite dimensional irreducible representations are given by its minimal (left) ideals $I\,\subset\,\Lambda(M)$, with $\Lambda(M)\,\vee\,I\,\subset\,I$.
The decomposition of this algebra into minimal ideals can be characterized by a spectral set $P_j$ of $\vee$-idempotents (\cite{graf}) in $(\Lambda(M), \vee)$, i.e. a set of elements satisfying:
\begin{enumerate}
\item $\sum_jP_j\,=\,1$
\item $P_j\vee P_k\,=\,\delta_{jk}P_j $
\item the rank of $P_j$ is minimal (non trivial), where the rank of $P_j$ is given by the dimension of the range of the  $\Lambda(M)$-morphism $\psi\,\mapsto\,\psi\vee P_j$.
\end{enumerate}
Under these conditions one has that 
\beq 
\label{dI}
I_j\,=\,\{\psi\,\in\,\Lambda(M)\,:\,\psi\vee P_j=\psi\}
\eeq
 and $\Lambda(M)\,=\,\oplus_jI_j$.
Elements in $I_j$ are called {\it spinors} with respect to the idempotent $P_j$. 
The action of the Dirac operator $\D$ 	defined in \eqref{difi} is meaningful on a set of spinors, i.e. $\D$ maps elements of $I_j$ into elements in $I_j$ with $I_{j}$ the left ideals of the Clifford algebra $(\Lambda(M), \vee)$, if and only if the condition
\beq
\label{compD}
P_j\vee\nabla_aP_j\,=\,0
\eeq
for any $\nabla_{a}\,=\,\nabla_{\del_a}$ holds.

The problem of defining a Clifford algebra and a  Hodge - de Rham Dirac operator on quantum spaces, i.e. for non commutative algebras deforming the algebras of functions on classical spaces,   has been widely studied following different approaches since for such non commutative spaces  there exists no canonical definition of  exterior algebra, differential calculus and symmetric tensors. 

Examples  of spin$^c$-spectral triples with a Dirac operator given by a non commutative deformation of the Dolbeault - Dirac operator (which classically coincides with the Hodge - de Rham operator on K\"ahler manifolds) has been  introduced and studied 
on quantum flag manifolds \cite{kraehmer}, quantised symmetric spaces \cite{kr-tu, matassa}:  via a suitable twist a spin$^c$-structure on quantum complex projective spaces $\C P_{q}^N$  is seen to give a spin-structure for odd $N$ in  \cite{dandrea-dabrowski, dandrea-landi}.

Within the quantum group formalism, meaningful Clifford algebras and spinors are introduced in \cite{frt, fiore} in terms of the properties of the $R$-braiding for the FRT  approach.  This approach is evolved in \cite{cridur1, cridur2, hecks}  for quantum groups equipped with a Woronowicz bicovariant exterior calculus, thus allowing for the definition of a Dirac operator.      
The papers \cite{heckschm, heckhodge, heckspin} develop a consistent formulation of the notions of Clifford algebras  for quantum groups equipped with left covariant Woronowicz calculi: the corresponding  exterior algebras are left modules for the Clifford algebra  (generalising  \eqref{claac}), spinors are introduced algebraically in terms of irreducibility subspaces of such an action. A consistent notion of metric tensors and Levi Civita connection acting  upon exterior forms is introduced,  a quantum analogue of  the Hodge - de Rham Dirac operator generalising \eqref{did} is defined and studied.  

The formalism  developed in this series of papers is indeed consistent under conditions which   are not satisfied for the case of the quantum group $\SU$ equipped with the family of three dimensional left covariant Woronowicz calculi studied in \cite{hec2001}. The aim of the present paper is to introduce a Hodge - de Rham Dirac operator on $\SU$ equipped with one example from that family of left invariant calculi. 

In section \ref{S3} we recall the classical construction for the Hodge - de Rham
operator on ${\rm S}^3$, and present its spectrum. In section \ref{sec:q}, instead of defining a Clifford algebra corresponding to a suitable symmetric bilinear form, we shall introduce a scalar product among forms, a Hodge duality operator for a meaningful quantum version of the classical Cartan-Killing metric tensor, and analyse a corresponding quantum deformation of the classical expressions \eqref{kmtco}, \eqref{diag}. We eventually write the explicit action of a quantum Hodge - de Rham operator $\D_{(q)}$, and compute its spectrum.

\section{The Hodge - de Rham operator  on ${\rm S^3}$}
\label{S3}

An explicit analysis of the Hodge - de Rham operator on ${\rm S}^3$ has been given in \cite{bangalore}. From that paper we recall the notations and the main results. 
We describe the sphere ${\rm S}^3$ as the ${\rm SU}(2)$ Lie group manifold, where points are described as the entries of the matrix 
\beq
\label{sug}
\gamma\,=\,\left(\begin{array}{cc} u & -\bar{v} \\ v & \bar{u}\end{array}\right), 
\qquad\bar{u}u+\bar{v}v=1.
\eeq
From a basis $T_a$ of the Lie algebra 
  $\mathfrak{su}(2)$, with
\beq
\label{cT}
[T_a, T_b]\,=\,\epsilon_{ab}^{\,\,\,\,c}T_c,
\eeq
the Maurer-Cartan 1-form
\beq
\label{dth}
\gamma^{-1}\dd \gamma\,=\,T_a\otimes\theta^a
\eeq
implicitly defines a left-invariant  basis $\{\theta^{a}\}_{a=1,\ldots,3}$ of $\Lambda^1({\rm SU}(2))$, whose dual vector fields are $\{L_a\}_{a=1,\ldots,3}$ representing the Lie algebra structure given by $[L_a,L_b]=\epsilon_{ab}^{\,\,\,\,c}L_c$. They have the following explicit expression:
\begin{align}
&\theta^x\,=\,i(u\dd v-v\dd u+\bar{v}\dd\bar{u} - \bar{u}\dd \bar{v}), \nn \\
&\theta^y\,=\,u\dd v-v\dd u-\bv\dd\bu+\bu\dd\bv, \nn \\
&\theta^z\,=\,2i(\bu\dd u+\bv\dd v)\,=\,-2i(u\dd\bu+v\dd\bv), \label{tfl}
\end{align}
and
\begin{align}
&L_x\,=\,\frac{i}{2}\,(\bv\del_u-v\del_{\bu}-\bu\del_v+u\del_{\bv}), \nn \\
&L_{y}\,=\,\frac{1}{2}\,(-\bv\del_u-v\del_{\bu}+\bu\del_v+u\del_{\bv}), \nn \\
&L_z\,=\,-\frac{i}{2}\,(u\del_u-\bu\del_{\bu}+v\del_v-\bv\del_{\bv}). \label{livf}
\end{align}
The ladder operators
\beq
L_{\pm}\,=\,\frac{1}{\sqrt{2}}\,(L_x\pm iL_y)
\label{lpm}
\eeq
with 
\beq
\label{lapm}
[L_-,L_+]=iL_z, \qquad[L_-,L_z]=-iL_-, \qquad[L_+,L_z]=iL_+,
\eeq
give an alternative left invariant basis, with  dual 1-forms given by
\beq
\theta^{\pm}\,=\,\frac{1}{\sqrt{2}}\,(\theta^x\mp i\theta^y).
\label{l1pm}
\eeq
The differential calculus is finally characterized by 
\beq
\dd\theta^a\,=\,-\frac{1}{2}f^{\,\,\,a}_{bc}\theta^b\wedge\theta^c
\label{mceq}
\eeq
in terms of the Lie algebra structure constants along a given left invariant basis,
with
\beq
\label{Lth}
L_a\theta^b\,=\,-\,f_{ac}^{\,\,\,\,\,b}\theta^c.
\eeq
The Cartan-Killing metric tensor
\beq
\label{ckme}
g\,=\,\theta^x\otimes\theta^x+\theta^y\otimes\theta^y+\theta^z\otimes\theta^z\,=\,\theta^-\otimes\theta^++\theta^+\otimes\theta^-\,+\,\theta^z\otimes\theta^z
\eeq
associated to the Lie algebra structure of $\mathfrak{su}(2)$ gives a riemannian metric tensor, whose corresponding \eqref{hod} basic 
Hodge duality map is, with $\tau\,=\,\theta^x\wedge\theta^y\wedge\theta^z$, 
\beq
\label{hos3}
\star 1\,=\,\tau, \qquad\qquad\star \theta^a\,=\,-\dd\theta^a
\eeq
with $\star^2\,=\,1$ on $\Lambda(S^3)$.


Upon the 4-dimensional  set $I_P$ of algebraic spinors with basis given by  (here $a\,=\,1, \ldots, 3$)
\begin{align}
&w_0\,=\,1-i\tau, \nn \\
& w_a\,=\,\delta_{ab}(\theta^b\,-\,\frac{i}{2}\,\varepsilon^b_{\,\,\,\,st}\theta^s\wedge\theta^t)\,=\,\delta_{ab}(\theta^b\,+\,i\dd\theta^b)\,=\,\delta_{ab}(\theta^b\,-\,i\star\theta^b),
\label{bases3}
\end{align}
the  action of the 
 Hodge - de Rham Dirac operator $\D\,=\theta^a\vee\nabla_a:\,I_P\,\to\,I_P$ has the following matrix form
\beq
\label{Dih}
\D\psi\,=\,\left(\begin{array}{cccc} 0 & L_+ & L_z & L_- \\ L_- & L_z-i & -L_- & 0 \\ L_z  & -L_+ & -i & L_- \\ L_+ & 0 & L_+ & -i-L_z \end{array}\right)\,\left(\begin{array}{c} \psi_0 \\ \psi_- \\ \psi_z \\ \psi_+ \end{array}\right)
\eeq
along the basis $w_0, w_z, w_{\pm}$. 
We consider the spectrum of the operator $\D$ acting upon square integrable spinors on $\rm{SU(2)}$, that is elements $\psi\in\C^4\otimes\mathcal{L}^2({\rm SU}(2)$. A Peter-Weyl basis for ${\rm SU}(2)$ is given by the elements 
\beq
\label{pwba}
f_m(J,N)\,=\,(R_+)^m\{v^{J-N/2}\bar{u}^{J+N/2}\},
\eeq
where $N\,\in\,\Z$, $J=s+|N|/2$ with $s\,\in\,\N$, and $m=1,\ldots, 2J+1$. This expression comes as the range of a  repeated action of the ladder operators $R_{\pm}$ (the right invariant vector fields on ${\rm SU}(2)$ corresponding to the elements $T_{\pm}$ in the Lie algebra \eqref{cT}) upon $f_0(J,N)=v^{J-N/2}\bar{u}^{J+N/2}$. For each fixed $J$, one has $\mathcal{L}^2({\rm SU}(2))\,\supset\,W_J=\{\phi\,: \,L^2\phi=-J(J+1)\phi\}$, with $\rm{dim}\,W_J=(2J+1)^2$. The spectrum of $\D$ is given by:
\begin{enumerate}
\item for $2J\neq\pm N$, one has 
\beq
\label{autov}
\lambda_{\pm}=\pm i\sqrt{J(J+1)}, \qquad\qquad 
\psi_{\pm}\,=\,\left(\begin{array}{c} \lambda_{\pm}\,f_m(J,N)  \\ L_-f_m(J,N) \\ L_zf_{m}(J,N) \\ L_+f_m(J,N) \end{array}\right),
\eeq
\beq
\label{eigtre}
\lambda=iJ, \qquad\qquad 
\psi\,=\,\left(\begin{array}{c} 0 \\ \frac{i}{J-N/2}\,L_{-}f_m(J,N) \\ f_m(J,N) \\ -\frac{i}{J+N/2}\,L_+f_m(J,N) \end{array}\right), 
\eeq
\beq
\label{eigv4}
\lambda=-i(J+1), \qquad\qquad 
\psi\,=\,\left(\begin{array}{c} 0 \\ -\frac{i}{1+J+N/2}\,L_{-}f_m(J,N) \\ f_m(J,N) \\ \frac{i}{1+J-N/2} \,L_+f_m(J,N) \end{array}\right); 
\eeq
\item for $2J=N$ one has $L_-f_m(J,2J)=0$, and 
\beq
\label{al1}
\lambda_{\pm}=\pm i\sqrt{J(J+1)}, \qquad\qquad\psi_{\pm}=\left(\begin{array}{c}\lambda_{\pm}f_m(J,2J) \\ 0 \\ L_zf_m(J,2J) \\ L_+f_m(J,2J)\end{array}\right),
\eeq
\beq
\label{al2}
\lambda=-i(J+1), \qquad\qquad\psi=\left(\begin{array}{c} 0 \\ 0 \\ -if_m(J,2J) \\ L_+f_m(J,2J)\end{array}\right),
\eeq
\beq
\label{al3}
\lambda=-i(J+1), \qquad\qquad\psi=\left(\begin{array}{c} 0 \\ 0 \\ 0 \\ f_m(J,2J)\end{array}\right);
\eeq
\item for $2J=-N$ one has $L_+f_m(J,-2J)=0$ and 
\beq
\label{al4}
\lambda_{\pm}=\pm i\sqrt{J(J+1)}, \qquad\qquad\psi_{\pm}\,=\,\left(\begin{array}{c}\lambda_{\pm}f_m(J,-2J)  \\ L_-f_m(J,-2J) \\ L_zf_m(J,-2J) \\ 0 \end{array}\right),
\eeq
\beq
\label{al5}
\lambda=-i(J+1), \qquad\qquad\psi\,=\,\left(\begin{array}{c} 0 \\ L_-f_m(J,-2J) \\ -i\,f_m(J,-2J) \\ 0 \end{array}\right),
\eeq
\beq
\label{al6}
\lambda=-i(J+1), \qquad\qquad\psi\,=\,\left(\begin{array}{c} 0 \\ f_m(J,-2J) \\ 0 \\ 0\end{array}\right).
\eeq
\end{enumerate}
By noticing that the action of $\D$ commutes with the action of $R_+$, we see that  $m$ is a degeneracy label for the eigenspinors corresponding to any given eigenvalue, so we have a basis for  $W_j\otimes \C^4$, i.e. a basis for   $I_P$ made by eigenspinors for $\D$.

\section{A Hodge - de Rham Dirac operator on quantum spheres}
\label{sec:q}

We move to the quantum setting  with a short presentation of  $\SU$ and a three dimensional left covariant differential calculus on it. This presentation is based on   \cite{wor87, ks}, more details can be found in  \cite{perugino}.  
As the quantum group $\SU$ we consider the  polynomial  unital  $*$-algebra $\ASU=(\SU,\Delta,S,\eps)$ generated by elements $a$ and $c$ such that, using the matrix notation
\beq
\label{Us}
 u = 
\left(
\begin{array}{cc} a & -qc^* \\ c & a^*
\end{array}\right) , 
\eeq
the Hopf algebra structure can be expressed as
\begin{align}
& uu^*=u^*u=1, \nn \\ &\Delta\, u = u \otimes u , \nn \\  &S(u) = u^* , \nn \\ &\eps(u) = 1. \label{hopfa} \end{align}
The deformation parameter
$q\in\IR$ is restricted  without loss of generality  to the interval $0<q<1$.  
By $\su$ we denote the universal envelopping algebra of $\ASU$, generated by the elements $E,F,K, K^{-1}$, with Hopf $^*$-algebra structure given by 
\begin{align}
 &K^{\pm}E=q^{\pm}EK^{\pm}, \nn \\
&K^{\pm}F=q^{\mp}FK^{\pm}, \nonumber \\
&   
[E,F] =\frac{K^{2}-K^{-2}}{q-q^{-1}},
\label{relsu}
\end{align} 
\beq
\label{stars}
K^*=K, \qquad  E^*=F ,
\eeq
and  
\begin{align}
&\Delta(K^{\pm}) =K^{\pm}\otimes K^{\pm}, \nn \\
&\Delta(E) =E\otimes K+K^{-1}\otimes E,  \nn \\ & 
\Delta(F)
=F\otimes K+K^{-1}\otimes F
\label{cosu}
\end{align}
for the coproduct, 
\begin{align}
& S(K) =K^{-1}, \nn \\
&S(E) =-qE, \nn \\
&S(F) =-q^{-1}F \label{ants}
\end{align}
for the antipode and  
\begin{align}
&\varepsilon(K)=1, \nn \\ 
& \varepsilon(E)=\varepsilon(F)=0
\label{cous}
\end{align}
for the counit. The algebra $\su$ is a Hopf $^*$-subalgebra included in the dual algebra $\ASU^o$, which is the largest Hopf $^*$-subalgebra contained in the dual vector space $\ASU^{\prime}$.  
The only non zero terms of the dual pairing are given upon the generators by 
\begin{align}
&K^{\pm}(a)=q^{\mp 1/2}, \qquad  K^{\pm}(a^*)=q^{\pm 1/2},  \nn \\ & E(c)=1, \qquad F(c^*)=-q^{-1}.\label{ndp}
\end{align}
Such a pairing gives rise to the 
 $^*$-compatible canonical commuting 
actions of $\su$ upon $\ASU$ 
\beq
h\lt x\,=\,x_{(1)}h(x_{(2)}),\qquad\qquad x\rt h=x_{(2)}h(x_{(1)}
\label{lra}
\eeq
for any $h\in\su$ and $x\in\ASU$.

A left covariant  first order differential calculi $(\Gamma, \dd)$ is a $\ASU={\mathcal H}$-bimodule provided 
$\dd:\ch\to\Gamma$ satisfies the Leibniz rule, $\dd (h\, h^{\prime}) = (\dd  h) h^{\prime}  + h\, \dd h^{\prime} $ for $h,h^{\prime}\in \ch$, and 
$\Gamma$ is generated by $\dd(\ch)$ as a $\ch$-bimodule. It is called  a $*$-calculus provided  there is an anti-linear involution $*:\Gamma\to\Gamma$ such that $(h_{1}(\dd h)h_{2})^*=h_{2}^{*}(\dd(h^*))h_{1}^{*}$ for any $h,h_{1},h_{2}\in \ch$. A first order differential calculus is said left covariant provided  a left coaction $\Delta_{L}^{(1)}:\Gamma\to\ch\otimes\Gamma$ exists, such that $\Delta_{L}^{(1)}(\dd h)=(1\otimes \dd)\Delta(h)$ and $\Delta_{L}^{(1)}(h_1\,\alpha\,h_2)=\Delta(h_{1})\Delta_{L}^{(1)}(\alpha)\Delta(h_2)$ for any $h,h_1,h_2\in\,\ch$ and $\alpha\in\,\Gamma$. The set $\Gamma$ turns out to be a free left covariant $\ch$-bimodule, with a free basis $\Gamma_{\rm L}$ of left invariant one forms,   namely the elements $\omega_{a}\in \Gamma$ such that $\Delta^{(1)}_{L}(\omega_{a})=1\otimes\omega_{a}.$ Its dimension is called the dimension of the first order calculus.
A first order  left covariant finite dimensional calculus is characterised by its quantum tangent space, namely a vector space ${\mathcal X}\subset\su$, which is dual to $\Gamma_{{\rm L}}$. 
A suitable basis of $\Gamma_{{\rm L}}$ gives a basis  of the quantum tangent space $\cx$ by duality, and  exact one-forms can be written as
\beq
\dd x\,=\,\sum_{a}(X_{a}\,\lt \,x)\omega_{a}
\label{ex1f}
\eeq
with $x\,\in\,\ASU$.

We consider in the present paper the three dimensional differential calculus on $\SU$ whose quantum tangent space $\cx$ has the basis
\begin{align}
\label{cx1}
&X_{z}\,=\,\frac{i}{2}\left(\frac{K^{-2}-1}{1-q}\right), \\ &X_{-}\,=\,-\frac{i}{\sqrt{2}}\,q^{-1/2}EK^{-1}, \nn \\ &X_{+}\,=\,-\frac{i}{\sqrt{2}}q^{1/2}FK^{-1}.
\end{align}
with corresponding  exact one forms
\begin{align*}
\dd a=-q\,c^*\omega_{+}\,+\,a\,\omega_{z},& \qquad\qquad \dd c=a^*\omega_{+}\,+\,c\,\omega_{z}, 
\\
\dd a^*=c\,\omega_-\,-\,q^{-1}a^*\omega_{z},& \qquad\qquad \dd c^*=-q^{-1}a\,\omega_{-}\,-\,q^{-1}c^*\omega_z.
\end{align*}
Such a differential calculus comes as a contraction of the well known $4D_+$
bicovariant differential calculus on $\SU$ introduced in \cite{wor89}. 
The basis is chosen so that antilinear hermitian conjugation upon them reads
$\omega_{-}^{*}\,=\,\omega_{+}, \; \omega_z^{*}\,=\,\omega_{z}$. 

Given the left covariant bimodule $\Gamma$ over $\ch$, an invertible linear mapping $\sigma:\Gamma\otimes_{\ch}\Gamma\,\to\,\Gamma\otimes_{\ch}\Gamma$ is called a braiding for $\Gamma$ provided $\sigma$ is a $\ch$-bimodule homomorphism which commutes with the left coaction on $\Gamma$ and satisfies the braid equation 
\beq
(1\otimes\sigma)\circ(\sigma\otimes 1)\circ(1\otimes \sigma)=(\sigma\otimes 1)\circ(1\otimes \sigma)\circ(\sigma\otimes 1)
\label{br3}
\eeq
on $\Gamma^{\otimes3}$. 
A braiding allows to define 
meaningful antisymmetriser operators $A^{(k)}:\Gamma^{\otimes k}\to\Gamma^{\otimes k}$. Under further suitable conditions their ranges suitably give an exterior differential algebra  $(\dd, \Gamma_{\sigma}, \wedge)$,  where the  exterior derivative $\dd$ is extended as a graded derivation with $\dd^2=0$,  satisfying a graded Leibniz rule (that is $\dd(\omega\wedge\omega^{\prime})=(\dd\omega)\wedge\omega^{\prime}+(-1)^{m}\omega\wedge\dd\omega^{\prime}$ for any $\omega\in\,\Gamma^{m}_{\sigma}$).

Such a braiding neither needs to exist nor it is unique for a given left covariant differential calculus over $\ch$: this is the main difference with bicovariant differential calculi, which present a  canonical braiding. 
A non canonical braiding for the calculus we are considering, whose range gives a well defined exterior algebra,  has been introduced in \cite{hec2001}. It is given by 
\begin{align}
\sigma(\omega_{\pm}\otimes\omega_{\pm})=\omega_{\pm}\otimes\omega_{\pm},& \nn  \\
\sigma(\omega_z\otimes\omega_z)=\omega_z\otimes\omega_z\,+\,\frac{2q(q-1)}{1+q^{-1}}(\omega_{+}\otimes\omega_-\,-\,\omega_{-}\otimes\omega_{+}),&
\nn \\
\sigma(\omega_+\otimes\omega_-)=q^2\omega_-\otimes\omega_+\,+\,(1-q^2)\omega_+\otimes\omega_-,&\qquad\qquad\sigma(\omega_{-}\otimes\omega_+)=\omega_+\otimes\omega_- \nn \\
\sigma(\oz\ot\om)=q^2\om\ot\oz\,+\,(1-q^2)\oz\ot\om,& \qquad\qquad \sigma(\om\ot\oz)=\oz\ot\om
\nn \\
\sigma(\op\ot\oz)=q^2\oz\ot\op\,+\,(1-q^2)\op\ot\oz,& \qquad\qquad \sigma(\oz\ot\op)=\op\ot\oz. \label{brad}
\end{align}
Such a braiding  has the following spectral decomposition,
\beq
(1-\sigma)(q^2+\sigma)=0,
\label{spbra} 
\eeq
with 
\begin{align}
&\dim\,\ker(1-\sigma)\,=6, \nn \\&\dim\,\ker(q^{2}+\sigma)\,=3.
\label{degen}
\end{align}
 The antisymmetriser operators  $A^{(k)}\,:\,\Gamma^{\otimes k}\to\Gamma^{\otimes k}$ they give rise  can be written as
\begin{align}
A^{(2)}=1-\sigma,\qquad\qquad&{\rm on}\,\Gamma^{\otimes 2} \nn \\ 
A^{(3)}=(1-\sigma_{2})(1-\sigma_1+\sigma_1\sigma_2),\qquad\qquad &{\rm on} \,\Gamma^{\otimes 3} 
\label{anespl}
\end{align}
with $\sigma_1=(\sigma\otimes 1)$ and $\sigma_2=(1\otimes \sigma)$, while 
$A^{(k)}$   are trivial for $k\geq4$. A basis of left invariant two forms in $\Gamma_{\sigma}^{2}$ is given by $\{\omega_{-}\wedge\omega_{+}, \omega_{+}\wedge\omega_{z}, \omega_{z}\wedge\omega_{-}\}$; 
given $\vartheta=i\omega_{-}\otimes\omega_{+}\otimes\omega_{z}$, left invariant volume forms are then, up to complex numbers,  
\beq
\theta\,=\,A^{(3)}(\vartheta).
\label{volf}
\eeq
The action of the antisymmetrisers $A^{(k)}$ on $\Gamma^{k}_{\sigma}$ is constant. Their  spectral resolution $A^{(k)}(\phi)\,=\,\lambda_{(k)}\phi$ with $\phi$ a $k$-form yields:
\begin{align}
&\lambda_{(2)}\,=\,(1+q^{2}) \nn \\ &\lambda_{(3)}\,=\,(1+q^2)(1+q^2+q^{4}). 
\label{spans}
\end{align}
The hermitian structure over left-invariant two forms and the wedge product antisymmetry read
 \begin{align*}
 (\omega_-\wedge\omega_+)^*&=-\omega_-\wedge\omega_{+}=q^{-2}\omega_{+}\wedge\omega_{-}, \\
 (\omega_-\wedge\omega_z)^*&=-\omega_z\wedge\omega_{+}=q^{-2}\omega_{+}\wedge\omega_{z}, \\
(\omega_+\wedge\omega_z)^*&=-\omega_z\wedge\omega_{-}=q^{2}\omega_{-}\wedge\omega_{z},
\end{align*}
while the volume form $\theta\in\Gamma^{3}_{\sigma}$ turns out to be a multiple of the classical one, namely the one we would obtain if  the braiding were the  classical flip,
\beq
\label{vol1}
\theta=q^2(\omega_-\otimes(\omega_{+}\otimes\omega_{z}\,-
\omega_z\otimes\omega_{+})\,+\,\op\otimes(\omega_z\otimes\omega_{-}\,-\,\om\ot\oz) 
\,+\,\oz\ot(\om\ot\op\,-\,\op\ot\om)).
\eeq
If a consistent braiding is defined on a left covariant first order differential calculus, then the duality $\Gamma_{{\rm L}}$ and $\cx$  allows to induce  a consistent
braiding $\sigma^{\dagger}$ on $\cx_{\cq}$. If 
\beq
\label{braind}
\sigma(\omega_a\otimes\omega_b)=\sigma_{ab}^{\,\,\,ks}\omega_k\otimes\omega_s
\eeq
then one has
\beq
\label{bratra}
\sigma^{\dagger}(X_a\otimes X_b)=\sigma_{ks}^{\,\,\,ab}X_k\otimes X_s,
\eeq
and one can define \cite{hec2001} a commutator $\left[~,~\right]:\cx_{\cq}^{\otimes 2}\to\cx_{\cq}$, fullfilling a peculiar Jacobi identity on the quantum tangent space corresponding to the calculus.  
For the calculus we are considering it is 
\begin{align}
&[X_a,X_a]=0, \nn \\
&[X_-,X_+]=-[X_+,X_-]=-i\left(\frac{2q^2}{1+q}\right)X_z, \nn \\
&[X_z,X_-]=-[X_-,X_z]=i\left(\frac{1+q}{2}\right)\,X_-, \nn \\
&[X_+,X_z]=-[X_z,X_+]=i\left(\frac{1+q}{2}\right)\,X_+.
\label{lieX}
\end{align}
If  the quantum commutator structure $\left[~,~\right]:\cx\times\cx\to\cx$ associated to the given braiding $\sigma$ 
is written as  $\left[X_a,X_b\right]=f_{ab}^cX_c$, one proves that a quantum version of the Maurer - Cartan formula for the differential calculus is valid, namely
\beq
\label{mcqu} 
\dd\omega_a=-\frac{1}{\lambda_{(2)}}f_{bc}^a\omega_b\wedge\omega_c.
\eeq
Together with the graded Leibniz rule, this relation characterises the action of the exterior differential $\dd$ upon the whole exterior algebra $\Gamma_{\sigma}$.

\subsection{A Hodge duality on $\SU$}
\label{ss:hodge}

As in \cite{perugino}, we introduce a Hodge duality operator on $\Gamma_{\sigma}$ in terms of suitable contractions, following the classical formulation that we recall. 

  Given a finite $N$-dimensional (complex, say) vector space $V$ the action of a bilinear form $\gamma:V\times V\to \C$ can be extended to give (for $j\geq s$) a map $\gamma:V^{\otimes s}\times V^{\otimes j}\to V^{\otimes j-s}$
by 
\beq
\label{cobi}
\gamma(v_{a_1}\otimes\ldots\otimes v_{a_s},v_{b_1}\otimes\ldots\otimes v_{b_j})=\gamma(v_{a_1},v_{b_1})\cdots \gamma(v_{a_s},v_{b_s})v_{b_{s+1}}\otimes\ldots\otimes v_{b_j}
\eeq
If a meaningful braiding $\sigma:V^{\otimes 2}\to V^{\otimes 2}$ is defined -- with $A^{(s)}:V^{\otimes s}\to V^{\otimes s}$ the corresponding antisymmetrisers -- then the bilinear form $\gamma$ allows to define a map over the corresponding exterior algebras, that is  $\gamma:V_{\sigma}^{s}\times V^{j}_{\sigma}\to V^{j-s}_{\sigma}$ via\footnote{We are denoting by $V_{\sigma}^s$ the range of $V^{\otimes s}$ under the action of the antisymmetriser $A^{(s)}$.} 
\beq
\label{exgam}
\gamma(v_{a_1}\wedge\ldots\wedge v_{a_{s}}, v_{b_1}\wedge\ldots\wedge v_{b_{j}})\,=\,\gamma(A^{(k)}(v_{a_1}\otimes\ldots\otimes v_{a_{k}}),A^{(s)}( v_{b_1}\otimes\ldots\otimes v_{b_{s}})).
\eeq
It is immediate to see that, if $\theta=v_1\wedge\ldots\wedge v_N$ denotes an element in the one dimensional vector space $V^N_{\sigma}$, upon defining $\tau=t\theta$ with $\C\ni t\neq 0$, for the map $\gamma_{\tau}:V^s_{\sigma}\to V^{N-s}_{\sigma}$ given by $\epsilon\mapsto\gamma_{\tau}(\epsilon)=(1/s!)\gamma(\epsilon, \tau)$ one has, with  $\epsilon\in V_{\sigma}^k$,  
\beq
\label{symho}
\gamma_{\tau}(\gamma_{\tau}(\epsilon))\,=\,t^2\,({\rm det}_{\theta}\gamma)(-1)^{k(N-k)}\epsilon\qquad\Leftrightarrow\qquad\gamma(v, v^{\prime})=\gamma(v^{\prime}, v) \quad\forall v,v^{\prime}\,\in\,V,
\eeq 
where one has defined 
$$
{\rm det}_{\theta}\gamma\,=\,\frac{1}{N!}\gamma(\theta, \theta):
$$
the degeneracy of the map $\gamma_{\tau}\circ\gamma_{\tau}$ is related to the symmetry of the bilinear $\gamma$.  Such a bilinear $\gamma$ turns out to be real, that is $\gamma(v,v^{\prime})\in\R$ for all $v^*=v, \,v^{\prime*}=v^{\prime}$, if and only if, with $\tau^*=\tau$, one has $\gamma_{\tau}(\epsilon^*)=(\gamma_{\tau}(\epsilon))^*$. it is evident that, upon a suitable choice for the global normalization factor $t$, the map $\gamma_{\tau}$ is a Hodge duality on $V_{\sigma}$ (the comparison with \eqref{hod} is straightforward). 

We extend this path to the quantum group setting, by considering the bilinear  map
$\gm\,:\,\Gamma_{\rm L}\times\Gamma_{\rm L}\,\to\,\IC$,  where $\Gamma_{\rm L}$ is the left invariant basis of 1-forms for the given exterior algebra $(\Gamma_{\sigma}, \wedge, \dd)$. 
Generalising \eqref{exgam}, this contraction is naturally extended to the left invariant part of $\Gamma_{\sigma}$ via $$
\gm(\omega_{a_1}\wedge\ldots\wedge\omega_{a_{k}}, \omega_{b_1}\wedge\ldots\wedge\omega_{b_{s}})\,=\,\gm(A^{(k)}(\omega_{a_1}\otimes\ldots\otimes\omega_{a_{k}}),A^{(s)}( \omega_{b_1}\otimes\ldots\otimes\omega_{b_{s}})). 
$$
From  \eqref{volf} 
we define the quantum volume form $\tau=\tau^*=\delta\theta$ with $\R\ni\delta\neq0$ and the determinant of the contraction $\gm$
\beq
\label{qdete}
{\rm det}_{\theta}\gm\,=\,\frac{1}{\lambda_{(3)}}\,\gm(\theta,\theta).
\eeq
We then define the linear operator  
$S:\Gamma^{k}_{\sigma}\to\Gamma^{3-k}_{\sigma}$ as  
\beq
\label{defS}
S(\omega)\,=\,\frac{1}{\lambda_{(k)}}\,\gm(\omega, \tau )
\eeq
on a left-invariant basis: the action of $S$ is extended upon all the exterior algebra $\Gamma^{\wedge}_{\sigma}$ by requiring the left $\ASU$-linearity. When $\gm$ is non degenerate and the eigenspaces for the action  of $S^2$ coincide with those of the action of the antisymmetrisers of the quantum differential calculus, then we say that the map $\gm$  is symmetric: such map is defined real if $S(\omega^*)=(S(\omega))^*$.
Under such symmetry conditions for $\gm$ we define $S$ to be a Hodge duality. 

 The Hodge duality allows to define \cite{kmt} an hermitian left invariant inner product on $\Gamma_{\sigma}$. This definition is based on the notion of integral, depending on a volume,  of an exterior form. One has 
$\int_{\tau}\omega=0$ if $\omega\in\Gamma^{k}_{\sigma}$ with $k\neq3$, while, for $\Gamma^{3}_{\sigma}\ni\omega=f\tau$ with $f\in\ASU$, one sets $\int_{\tau}\omega=h(f)$
with $h:\ASU\to\IC$ the Haar state defined on $\SU$.  In terms of this definition of integral, one sets
\beq
\label{digr}
\hs{\omega}{\omega^{\prime}}=0 
\eeq
when $\omega$ and $\omega^{\prime}$ have different degrees, while, for $\omega, \omega^{\prime}$ both in $\Gamma^{ k}_{\sigma}$, one defines
\beq
\label{eqdi}
\hs{\omega}{\omega^{\prime}}=\int_{\tau}\omega^{\prime*}\wedge S(\omega).
\eeq
The inner product defined above is left invariant: one can see that 
\beq
\label{liip}
\hs{f\omega}{f^{\prime}\omega^{\prime}}=h(f^{\prime*}f)\hs{\omega}{\omega^{\prime}}
\eeq
for any $f,f^{\prime}$ in $\ASU$ and left invariant exterior forms $\omega, \omega^{\prime}$.

Given the inner product \eqref{eqdi} it is possible to define an adjoint operator to the exterior differential $\dd$, namely $\dd^{\dagger}:\Gamma^{ k}_{\sigma}\to\Gamma^{ k-1}_{\sigma}$ satisfying the relation 
\beq
\label{add}
\hs{\dd^{\dagger}\omega}{\omega^{\prime}}=\hs{\omega}{\dd \omega^{\prime}}
\eeq
for any $\omega, \omega^{\prime} \in\,\Gamma_{\sigma}$. It turns out \cite{kmt} that the action 
of the operator $\dd^{\dagger}$ is
\beq
\label{diago}
\dd^{\dagger}\omega=(-1)^kS^{-1}\dd (S\omega)
\eeq
for $\omega\in\Gamma^{ k}_{\sigma}$. Notice that, by identifying the operator $S$ to a Hodge duality, the expression for the action of $\dd^{\dagger}$ in the quantum setting coincides with that in the classical setting, see \eqref{kmtco}, \eqref{diag}. Notice as well that the expression \eqref{diago} holds for any $S$ corresponding to a non degenerate $\gm$, irrespective of the degeneracy of the spectrum  of the operator $S^2$.

\subsection{A Cartan - Killing metric on $\SU$}
\label{cksuq}

Among all possible bilinear maps  we select those given by (non degeneracy being equivalent to $\alpha\,\beta\,\gamma\,\neq\,0$)
\begin{align}
\label{g1fo}
&\gm(\omega_{-},\omega_{+})=\alpha, \\ &\gm(\omega_{+},\omega_{-})=\beta, \nn \\ & \gm(\omega_{z},\omega_{z})=\gamma.
\end{align}
From \cite{perugino} we know that the map $S$ is well defined for all the three dimensional left invariant calculi $\mathcal{K}$ on $\SU$ studied in that paper,  the calculus we are considering in the present paper being one of them. In particular, the condition that  the degeneracy of the spectrum of $S^2$ is maximal and that the action of $S$ commutes with the hermitian conjugation amount to linearly constrain
 the real variables $\alpha$ and $\beta$. Following \cite{warped}, we set the condition (with $\IR\ni\xi\neq0$)
 \beq
 \label{qck}
S(\omega^a)\,=\,\xi\dd\omega^a
\eeq
for any left invariant one form $\omega^a$ -- which mimics the classical condition \eqref{hos3} that characterises the Cartan - Killing metric tensor up to a global scaling -- as the definition of the quantum version on $\SU$ of the Cartan - Killing metric. 

From \eqref{vol1} the 
 one has   
$\gm(\theta,\theta)=\,6q^4(\alpha\beta\gamma)$.
The action of the operator $S:\Gamma^{k}_{\sigma}\to\Gamma^{ 3-k}_{\sigma}$ defined in \eqref{defS} is 
\begin{align}
S(\omega_-)=i\delta\alpha\,\omega_z\wedge\omega_-, & \qquad S(\omega_z\wedge\omega_-)=
-2iq^4(\delta/\lambda_{(2)})\alpha\gamma\,\omega_- \nn \\
S(\omega_+)=i\delta\beta\,\omega_+\wedge\omega_z, & \qquad S(\omega_+\wedge\omega_z)=
-2iq^4(\delta/\lambda_{(2)})\beta\gamma\,\omega_+ \nn \\
S(\omega_z)=i\delta q^2\gamma\,\omega_-\wedge\omega_+, & \qquad S(\omega_-\wedge\omega_+)=
-2iq^2(\delta/\lambda_{(2)})\alpha\beta\,\omega_z 
\label{Sge}
\end{align}
and also
\begin{align}
&S(1)\,=\,\tau, \nn \\
&S(\tau)\,=\,\frac{6q^4}{\lambda_{(3)}}\delta^2\,\alpha\beta\gamma
\label{sut}
\end{align}
From the relations above it is immediate to see that $S^2(\omega_a)=\zeta\omega_a$ for any left invariant 1-form if and only if $\alpha=\beta$, and that $(S(\omega))^*=S(\omega^*)$ if and only if $\alpha, \beta, \gamma$ are real parameters. One sees that the Cartan - Killing metric ${\rm g}$ (among those defined in \eqref{g1fo})  set by \eqref{qck} is characterised by 
\beq
\label{cqk}
\alpha=\beta=\left(\frac{1+q}{2q^2}\right)^2\gamma.
\eeq
We shall refer to the $S$ operator defined by \eqref{defS}, \eqref{g1fo}, \eqref{cqk} with $\alpha\in\IR$ as the Hodge duality operator $\star:\Gamma^{k}_{\sigma}\to\Gamma_{\sigma}^{3-k}$ corresponding to the \emph{quantum} Cartan - Killing metric. We select the normalization in such a way that $S(\tau)=1$, so that 
\begin{align}
&\star(1)=\tau, \nn \\
&\star(\tau)=1, \nn \\
&\star^2(\omega)=A\omega, \nn \\
&\star(\omega_a)=\xi\dd\omega_a
\label{caho}
\end{align}
for any 1- and 2-form $\omega$, and any left invariant 1-form $\omega_a$. Notice that, given the relations \eqref{cqk}, one has 
\begin{align}
&A=\left(\frac{2}{\lambda_{(2)}}\right)\left(\frac{\lambda_{(3)}}{6}\right), \nn \\
&\xi=-q^{-2}\left(\frac{\lambda_{(3)}}{6}\right)^{1/2}\frac{1}{\sqrt{\alpha}}.
\label{coef}
\end{align}
The coefficient $A$ is related to the deformation of the spectrum of the quantum braiding with respect to the classical flip. 

\subsection{A Hodge - de Rham operator on $\SU$}
\label{sec:hdqu}

The aim of this section is to prove that a suitable quantum deformation $\D_{(q)}$ of the action of the Hodge - de Rham operator \eqref{difi} can be completely reduced on a subspace in $\Gamma_{\sigma}$
spanned by a suitable deformation of the elements \eqref{bases3}. We define the linear operator 
\beq
\Gamma^{ k}_{\sigma}\ni\omega\quad\mapsto\quad\D_{(q)}(\omega)\,=\,\dd\omega+\varepsilon(k)\star\dd(\star\omega)
\label{dedD}
\eeq
where $\epsilon(k)$ is a complex parameter depending on the degree $k=1,\ldots,3$.
We consider the element
\beq
\label{psize}
\Gamma_{\sigma}\ni\psi_0=1+\kappa\tau:
\eeq
it is immediate to see that the relation, with $f\in\ASU$, 
\beq
\label{Dqz}
\D_{(q)}:f\psi_0\quad\mapsto\quad(X_a\lt f)\{\omega_a+\varepsilon(3)\kappa(\star\omega_a)\}
\eeq
holds. From this result we define the three elements
\beq
\Gamma_{\sigma}\ni\phi_a=\omega_a+\varepsilon(3)\kappa(\star\omega_a)
\label{basq}
\eeq
and compute the expression $\D_{(q)}(f_a\phi_a)$, namely the action of $\D_{(q)}$ upon elements in the subspace  (a left $\ASU$- module in $\Gamma_{\sigma}$) spanned by \eqref{basq}. If we denote by $I\subset\Gamma_{\sigma}$ the free left $\ASU$-module spanned by $\{\psi_0, \phi_a\}$, it is easy, using \eqref{caho}, to prove that $\D_{(q)}$ maps elements in $I$ to elements in $I$ if and only if the conditions
\begin{align}
&\varepsilon(1)=\varepsilon(3), \nn \\
&\varepsilon(2)(\varepsilon(3))^2\kappa^2A^2=1
\label{conpro}
\end{align}
hold. 
Since only the $A$ parameter is fixed in the above relations, the space of solutions of such constraints is quite rich. If one wants to mimic as much as possible the classical expression, then one can choose
\beq
\label{soluz} 
\varepsilon(1)=\varepsilon(3)=-\varepsilon(2)=1, \qquad\kappa=\pm iA^{-1}.
\eeq
We choose $\kappa=-iA^{-1}$. The Dirac operator \eqref{dedD} under the conditions \eqref{soluz} has a matrix representation along the basis $\{\psi_0, \phi_-,\phi_z,\phi_+\}$ on the space $I$ given by (the action $X\lt f$ is defined in \eqref{lra}, with  $\alpha, \beta, \gamma$ satisfying the constraints \eqref{cqk})
\beq
\label{dirq}
\D_{(q)}\,\left(\begin{array}{c} f_0 \\ f_- \\ f_z \\ f_+\end{array}\right) \,=\, \left(\begin{array}{cccc} 
0 & q^4\alpha X_+ & q^2\gamma X_z  &  \beta X_- \\ X_- & A\xi^{-1}(i-\frac{2q^2}{1+q}X_z) & A\xi^{-1} \frac{2}{1+q}X_- & 0 \\ X_z & A\xi^{-1} \frac{(1+q)}{2}X_+ & A\xi^{-1} (i+(1-q)X_z) & -A\xi^{-1}\frac{1+q}{2q^2}X_- \\  X_+ & 0 & -A\xi^{-1}\frac{2q^2}{1+q}X_+ & A\xi^{-1} (i+\frac{2}{1+q}X_z) 
\end{array} \right)\lt
\left(\begin{array}{c} f_0 \\ f_- \\ f_z \\ f_+\end{array}\right)
\eeq
The limit $q\to1$ the action of the operator $\D_{(q)}$ coincides\footnote{Notice that, in the classical limit $q\to1$, one sees from \eqref{coef} that $A\to1$ and $\xi\to-1$.} with the action of the operator \eqref{Dih}. 
Our aim is now to determine the spectrum of the operator $\D_{(q)}$ when acting upon $I$, that we realise as the quantum analogue of the space $\C^4\otimes{\mathcal L}^2({\rm SU}(2))$ described in the section \ref{S3}, more precisely as the complex linear span along the quantum analogue of the classical Wigner functions for ${\rm SU}(2)$, as showed by the Peter-Weyl theorem. We realise such a basis, related to the so called PBW construction, as 
\beq
\label{qwig}
f_{m}(J,N)\,=\,(c^{J-N/2}a^{\star J+N/2})\rt E^m
\eeq
with $J=1/2,1,3/2,\ldots$, $N=-2J, -2J+1, \ldots, 2J-1, 2J$ and  $m=0, 1, \ldots, 2J$. It is evident that the existence of such a basis directly follows from the analogy between the theory of unitary representations of the Lie group ${\rm SU}(2)$ and the theory of unitary corepresentations \cite{ks} of the quantum group ${\rm SU}_q(2)$.

\begin{itemize}
\item 
In order to determine the spectrum of this operator we start by considering an \emph{Ansatz} close to the one used in \eqref{autov} (see \cite{bangalore}). We consider the spinor 
\beq
\label{trsp}
\psi\,=\,\left(\begin{array}{c} \sigma f_m(J, N) \\ \mu X_-\lt f_m(J,N) \\ \rho f_m(J,N) \\ \tilde{\mu} X_+\lt f_m(J,N)\end{array}\right).
\eeq  
The eigenvalue equation $\D_{(q)}\psi=\lambda\psi$  upon the spinor \eqref{trsp} is equivalent to the following algebraic matrix eigenvalue equation
\begin{align}
\label{eieq}
&\left(\begin{array}{cccc} 0 & -\frac{\alpha}{2}q^{3-N}\phi(J,N) & \frac{2iq^5}{1+q}\alpha\left[\frac{N}{2}\right]q^{-N/2} & -\frac{\alpha}{2}q^{1-N}\epsilon(J,N) \\  1 & -A\xi^{-1}iq^{-1}\left[\frac{N}{2}\right]q^{-N/2} & A\xi^{-1}\frac{2}{1+q} & 0 \\ 
\frac{i}{2}\left[\frac{N}{2}\right](1+q^{-1})q^{-N/2} & -A\xi^{-1}\frac{1+q}{4}q^{-(1+N)}\phi(J,N) & \frac{i}{2}A\xi^{-1}(1+q^{-N}) & A\xi^{-1}\frac{1+q}{4}q^{-(1+N)}\epsilon(J,N) \\ 1 & 0 & -A\xi^{-1}\frac{2q^2}{1+q} & A\xi^{-1}iq\left[\frac{N}{2}\right]q^{-N/2} 
\end{array}
\right)\cdot\, \\  & \qquad\qquad \qquad\qquad\qquad\qquad\qquad\qquad\qquad \qquad\qquad\qquad\qquad\qquad\cdot\left(\begin{array}{c} \sigma \\ \mu \\ \rho \\ \tilde{\mu} \end{array}\right)\, =\,\lambda
\,\left(\begin{array}{c} \sigma \\ \mu \\ \rho \\ \tilde{\mu} \end{array}\right),
\nn \end{align}
where
\begin{align}
\label{phia}
&EF\lt f_m(J,N)\,=\,\epsilon(J,N)\,f_m(J,N)\,=\,\left([J+\frac{1}{2}]^2-[\frac{1-N}{2}]^2\right)\,f_m(J,N) , \\
&FE\lt f_m(J,N)\,=\,\phi(J,N)\,f_m(J,N)\,=\,\left([J+\frac{1}{2}]^2-[\frac{N+1}{2}]^2\right)\,f_m(J,N) \nn 
\end{align} 
\beq
[x]\,=\,\frac{q^x-q^{-x}}{q-q^{-1}}.
\label{qn}
\eeq
The solutions of the eigenvalue problem for \eqref{eieq} are
\begin{align}
&\lambda_{\pm}=\pm i\left(\frac{\alpha}{2}\,q^{1-N}\left\{q^2\phi(J,N)+\epsilon(J,N)+2q^3[\frac{N}{2}]^2\right\}\right)^{1/2}, \nn \\
&\tilde{\mu}=\mu, \qquad \rho=\,\frac{i}{2}\,\mu\,[\frac{N}{2}](1+q^{-1})q^{-N/2} , \qquad \sigma=\,\mu\lambda_{\pm}, 
\label{qeig}
\end{align}
and 
\begin{align}
&\lambda^{\prime}_{\pm} = \frac{i}{2}\left( A\xi^{-1} \right)\left[ \frac{3-q^{-N}}{2} \pm \sqrt{\frac{(3-q^{-N})^2}{4}-2+2q^{-N}\left( \left[ J+1 \right]^2 + \left[ J \right]^2 \right)} \right], \nn \\
&\frac{\mu}{\rho} = \dfrac{\left( \lambda_{\pm} - i A\xi^{-1} q^{1-N/2}\left[ \frac{N}{2} \right] \right)}{\left( \lambda_{\pm} + i A\xi^{-1} q^{-1-N/2}\left[ \frac{N}{2} \right] \right)} \frac{\tilde{\mu}}{\rho}+ \dfrac{2A\xi^{-1}(1+q^2)}{\left( 1+q \right) \left( \lambda_{\pm} + i A\xi^{-1} q^{-1-N/2}\left[ \frac{N}{2} \right] \right)}, \nn \\
&\frac{\sigma}{\rho} = \left( \lambda_{\pm} - i A\xi^{-1} q^{1-N/2}\left[ \frac{N}{2} \right] \right) \frac{\tilde{\mu}}{\rho} + \dfrac{2 A \xi^{-1}q^2}{1+q}, \nn \\
&\frac{\tilde{\mu}}{\rho} = \{2\lambda^2 A\xi^{-1}\frac{q^2}{1+q}+ \frac{2iq^{1-N/2}}{1+q}\left[ \frac{N}{2} \right]((A\xi^{-1})^2-q^4\alpha)\lambda + \alpha \frac{q^{3-N}(A\xi^{-1})}{1+q} ( 1+q^2 + 2q\left[ \frac{N}{2} \right] ) \}\cdot
\nn \\
&
\qquad\qquad \left\{-\lambda_{\pm} ^3 -i \lambda_{\pm} ^2 A\xi^{-1}\left[ \frac{N}{2} \right] q^{-N/2}(q-q^{-1})- \lambda_{\pm} ( \frac{\alpha}{2} (q^{3-N} \phi + q^{1-N} \epsilon ) + (A\xi^{-1})^2 q^{-N}\left[ \frac{N}{2} \right]^2  ) + \right. \nn \\ &\qquad\qquad\qquad\left. i A\xi^{-1} \frac{\alpha}{2} q^{-3N/2}\left[ \frac{N}{2} \right] (q^4 \phi - \epsilon) \right\}^{-1}. \nn 
\end{align}
\item For $2J=N$ one has $X_-\lt f_m(J=N/2, N)=0$. We explore the eigenvalue problem upon considering the ansatz
\beq
\label{anjn}
\psi\,=\,\left(\begin{array}{c} \sigma f_m(J=N/2, N) \\ 0 \\ \rho f_m(J=N/2, N) \\ \tilde{\mu} X_+\lt f_m(J=N/2,N) \end{array}\right)
\eeq
for an eigenspinor of $\D_{(q)}$. The action of $\D_{(q)}$ is meaningful upon such a space of spinors, and the eigenvalue problem for $\D_{(q)}$ is equivalent to the following matrix eigenvalue problem 
\beq
\label{m3au}
\left(\begin{array}{ccc} 0 & \frac{i}{2}\,q^2\gamma\,\frac{q^{-N}-1}{1-q} & -\frac{1}{2}\,\beta q^{1-N}[N] \\ \frac{i}{2}\left(\frac{q^{-N}-1}{1-q}\right) & A\xi^{-1}\frac{i}{2}(1+q^{-N}) & A\xi^{-1}\frac{1}{4}\,q^{-1-N}(1+q)[N] \\ 1 & -A\xi^{-1}\,\frac{2q^2}{1+q} & iA\xi^{-1}q^2 \left(\frac{q^{-N}-1}{1-q^2}\right)
\end{array}\right)
\,\left(\begin{array}{c} \sigma \\ \rho \\ \tilde{\mu}\end{array}\right)\,=\,\lambda\,
\left(\begin{array}{c} \sigma \\ \rho \\ \tilde{\mu}\end{array}\right)
\eeq
The first set of eigenvalues and corresponding eigenspinors is given by 
\begin{align}
&\lambda_{\pm}=\pm i\left(\frac{\alpha}{2}\,q^{1-N}\left\{q^2\phi(J=N/2,N)+\epsilon(J=N/2,N)+2q^3[\frac{N}{2}]^2\right\}\right)^{1/2}\, \nn \\ &\qquad\qquad =\,\pm i\left(\frac{\alpha}{2}\,q^{1-2J}\left\{[2J]+2q^3[J]^2\right\}\right)^{1/2}, \nn \\
&\tilde{\mu}\neq0, \qquad \rho=\,\frac{i}{2}\,\mu\,[\frac{N}{2}](1+q^{-1})q^{-N/2} , \qquad \sigma=\,\tilde{\mu}\lambda_{\pm}.
\label{qeig2}
\end{align}
Such eigenvalues coincide, for $N=2J$, with the eigenvalues \eqref{qeig}: the parameter $\rho$ turns indeed out to be given by  $X_z\lt a^{\star N}=\rho\,a^{\star N}$. This allows to write the eigenspinors corresponding to the eigenvalues \eqref{qeig2} as 
\beq
\label{ausb}
\psi\,=\,\left(\begin{array}{c} \lambda_{\pm}f \\ 0 \\ X_z\lt f \\ X_+\lt f\end{array}\right)
\eeq  
con $f=f_m(J=N/2, N)$.

The third eigenvalue of the matrix \eqref{m3au} is given by 
\beq
\label{tauv}
\lambda=iA\xi^{-1}\left(\frac{1}{2}(1+q^{-2J})+q^2\left(\frac{q^{-2J}-1}{1-q^2}\right)\right),
\eeq
its corresponding eigenspinor is given by
\begin{align}
\label{eispqq}
& \sigma = \frac{i \sqrt{\frac{3}{2}} \sqrt{\alpha} \left(q^2+1\right) \sqrt{\frac{q^4+q^2+1}{q^2+1}} q^{-N} \left(q^{2 N}-1\right)}{5 q^{N+2}+3 q^{N+4}+4 q^N-q^4-3 q^2-2} \\
& \rho = -\frac{i (q+1) q^{-N-4} \left(3 q^{2 N}+2 q^{N+2}+2 q^{N+4}+2 q^{N+6}+7 q^{2 N+2}+5 q^{2 N+4}+3 q^{2 N+6}-q^6-3 q^4-5 q^2-3\right)}{4 \left(5 q^{N+2}+3 q^{N+4}+4 q^N-q^4-3 q^2-2\right)} \nn \\
&\tilde{\mu} = 1
\nn \end{align}
Once fixed $2J=N$, a direct inspection shows that the action of the operator $\D_{(q)}$ 
is diagonal upon the pair eigenvalue-eigenspinor given by  
\beq
\label{al3q}
\lambda=iA\xi^{-1}q^{-J}[1+J], \qquad\qquad\psi=\left(\begin{array}{c} 0 \\0 \\ 0 \\ f_m(J=N/2, N)\end{array}\right)
\eeq
\item We now consider the case $N=-2J$, with $X_+\lt f_m(J=-N/2, N)=0$, and consider for an eigenspinor the ansatz
\beq
\label{anjm}
\psi\,=\,\left(\begin{array}{c} \sigma f_m(J=-N/2, N) \\ \mu X_-\lt f_m(J=-N/2, N) \\ \rho f_m(J=-N/2, N) \\ 0 \end{array}\right).
\eeq
The eigenvalue problem for the action of $\D_{(q)}$ is meaningful upon such a space of spinors, and turns out to be equivalent to the matrix egenvalue problem given by 
\beq
\label{m3au2}
\left(\begin{array}{ccc} 0 & -\frac{1}{2}\alpha q^{3+2J}[2J] & q^2\gamma\,\frac{i}{2}\left(\frac{q^{2J}-1}{1-q} \right)\\ 1 & iA\xi^{-1}\left(\frac{1-q^{2J}}{1-q^2}\right) & A\xi^{-1}\frac{2}{1+q} \\ \frac{i}{2}\left(\frac{q^{2J}-1}{1-q}\right) & -\frac{1}{4}A\xi^{-1}(1+q)q^{{2J}-1}[2J] & A\xi^{-1}\frac{i}{2}(1+q^{2J}) 
\end{array}\right)
\,\left(\begin{array}{c} \sigma \\ \mu \\ \rho \end{array}\right)\,=\,\lambda\,
\left(\begin{array}{c} \sigma \\ \mu \\ \rho \end{array}\right)
\eeq
The first set of eigenvalues and corresponding eigenvectors for this matrix is
\begin{align}
&\lambda_{\pm}=\pm i\left(\frac{\alpha}{2}\,q^{1+2J}\left\{q^2[2J]+2q^3[J]^2\right\}\right)^{1/2}, \nn \\
&\mu\neq0, \qquad \rho=\,-\frac{i}{2}\,\mu\,[J](1+q^{-1})q^{J} , \qquad \sigma=\,\mu\lambda_{\pm}.
\label{qeig3}
\end{align}
The third eigenvalue of \eqref{m3au2} is 
\beq
\label{tera}
\lambda=iA\xi^{-1}\left(\frac{1}{2}(1+q^{2J})+\left(\frac{1-q^{2J}}{1-q^2}\right)\right),
\eeq
with corresponding eigenspinor 
\begin{align}
\label{espb}
&\sigma = \frac{2 \sqrt{6} \sqrt{\alpha} q^2 \left(q^2+1\right) \sqrt{\frac{q^4+q^2+1}{q^2+1}} \left(q^{2J}-1\right)}{(q+1) \left(3 q^{2J+2}+q^{2J+4}+2 q^{2J}+q^4-q^2\right)} \\
&\mu = -\frac{4 i \left(q^{2J+4}-q^{2J}+q^4+2 q^2+3\right)}{(q+1) \left(2 q^{4J}+2 q^{2J+2}+2 q^{2J+4}+3 q^{4J+2}+q^{4J+4}+2 q^{2J}+q^4-q^2\right)} \nn \\
&\rho = 1.
\nn
\end{align}
Analogously to \eqref{al3q}, the action of $\D_{(q)}$ is diagonal upon the pair eigenvalue-eigenspinor given by 
\beq
\label{al6q}
\lambda=iA\xi^{-1}\{1+q^{J+1}[J]\}, \qquad\qquad\psi=\left(\begin{array}{c} 0 \\ f_m(J=-N/2, N)\\ 0 \\ 0\end{array}\right)
\eeq
\end{itemize}
The spectrum of $\D_{(q)}$ turns out, by analysing the classical limit $q\to 1$, to be a quantum deformation of the spectrum of the Hodge - de Rham operator studied in \cite{bangalore} and described in the section \ref{S3}. Such a spectrum depends not only on $J$, as in the classical setting, but also on $N$. This is a well known phenomenon in the theory of differential operators on quantum spaces: quantising a classical differential operator amounts to remove some of the degeneracies of its classical spectrum. 

Upon explicitly counting the multiplicities of the the eigenvalues written above (notice that, as obvious since the action of $\D_{(q)}$ is written in terms of only left acting operators, its spectrum does not depend on $m$), one can, as already pointed out in section \ref{S3}  for the classical counterpart of $\D_{(q)}$, say that we have explicitly written a basis for $I$ made of eigenspinors for $\D_{(q)}$. 

We also mention that we equipped the quantum group $\SU$ with a well known left covariant three dimensional calculus, so that all the computations related to the explicit action of the Dirac operator and its spectrum strongly depends on it. We leave to a further analysis to present Hodge - de Rham Dirac operators on the same quantum group $\SU$ equipped with the family of three dimensional left covariant calculi studied in \cite{perugino}.

\subsection*{Acknowledgements} We should like to thank Patrizia Vitale and  Franco Ventriglia for many interesting discussions. We  gratefully acknowledge the support of the INFN, of the Spanish MINECO grant MTM2014-54692-P and Quitemad+, S2013/ICE-2801, of the University of Luxembourg. We worked at this paper while visiting the ICMAT (Madrid), the Center for High Energy Physics (Bangalore), the Brookhaven National Laboratory (Upton, NY): we are happy to thank our hosts there.

\end{document}